\documentclass[12pt]{amsart}      
\usepackage{amssymb,enumerate}
\usepackage{amsmath}
\usepackage{setspace}
\usepackage{tcolorbox}

\newtheorem{thm}{Theorem}[section] 
\newtheorem{lemma}[thm]{Lemma} 
\newtheorem{prop}[thm]{Proposition} 
\newtheorem{cor}[thm]{Corollary} 
\newtheorem{example}[thm]{Example} 
\theoremstyle{definition} 
\newtheorem{remark}[thm]{Remark} 
\newtheorem{remarks}[thm]{Remarks} 
\newtheorem{ques}[thm]{Question}
\newtheorem{definition}[thm]{Definition} 
\newtheorem{example2}[thm]{Example} 
 
\def\F{\mathcal F} 
\def\U{\mathcal U} 
\def\V{\mathcal V}

\def\diam{\operatorname{diam}} 
\def\co{\operatorname{co}} 
\def\en{\mathbb N} 
\def\er{\mathbb R} 
\def\qe{\mathbb Q}
\def\r{|} 
\def\ep{\varepsilon} 
\def \Bos {\operatorname{Bos}} 
 
\def\Hs{\operatorname{Hs}}
\def \Bof {\operatorname{Bof}} 
 
\def\Hf{\operatorname{Hf}}

\def\clust{\operatorname{clust}}
\def\ck{\operatorname{ck}} 
\def\sfrag{\operatorname{\sigma-frag}} 
\def\sfra{\operatorname{Fr}_\sigma} 
\def\frag{\operatorname{frag}}
 
\def\orank{\operatorname{\beta}} 
\def\card{\operatorname{card}} 
\def\Frag{\operatorname{Frag}}
\def\oorank{\operatorname{\omega_1-\beta}}
\def\range{\operatorname{range}}

\begin{document}

\title[Distances to spaces of first resolvable class mappings]{Distances to spaces of first resolvable class mappings}
\author[P. Ludv\'\i k]{Pavel Ludv\'\i k}

\address{Department of Mathematical Analysis and Applications of Mathematics\\
Faculty of Science Palack\'y University Olomouc \\
17. listopadu 12, 771 \ 46\\
Olomouc, Czech Republic}
\email{pavel.ludvik@upol.cz}

\subjclass[2020]{54E52, 54H05, 54C35}

\keywords{Baire and Borel functions, fragmentability, resolvable set, mappings of first $H$-class, oscillation rank}

\thanks{}

\begin{abstract}
We study the mappings of the first resolvable class defined by G. Koumoullis in \cite{koumoullis} as a valuable tool to address the point of continuity property in the non-metrizable setting. First, we investigate the distance of a general mapping to the family of mappings of the first resolvable class via the \emph{fragmentability} quantity. We partially generalize papers of B. Cascales, W. Marciszewski, M. Raja; C. Angosto, B. Cascales, I. Namioka; and J. Spurn\'{y}. Second, we introduce the class of mappings with the countable oscillation rank, study its basic properties and relate it to the mappings of the first resolvable class and other well known classes of mappings. This rank has been in a less general context considered by S.~A. Argyros, R. Haydon and some others.
\end{abstract}
\maketitle

\section{Introduction}
We open the section with an important notion of fragmented mapping introduced in \cite{koumoullis}.
Let $X$ be a topological space and $E$ a metric space. Then a mapping $f:X\to E$ is called \emph{$\ep$-fragmented}  if for every nonempty (equivalently, nonempty closed) set $F\subset X$ there exists an open set $V\subset X$ such that $V\cap F \neq \emptyset$ and $\diam f(V\cap F) \leq \ep$. We write briefly
\[
\frag (f) = \inf \{\ep>0: f \text{ is $\ep$-fragmented} \}
\]
and term $f$  to be \emph{fragmented} if $\frag(f)=0$. We designate the set of all fragmented mapping $f:X\to E$ by $\Frag(X,E)$.

A set $H$ in $X$ is \emph{resolvable} if its characteristic function $\chi_H$ is fragmented, that is, for any nonempty (equivalently nonempty closed) $A\subset X$ there exists a nonempty relatively open set $U\subset A$ such that either $U\subset H$ or $U\subset X\setminus H$. We refer the reader to \cite[\S12,V-VI]{kuratowski-top} or \cite[Section~2]{koumoullis} for elementary properties of resolvable sets. We just recall that the family $\Hs(X)$ of all resolvable sets forms an algebra containing all open sets.

Having a metric space $E$, the function $\diam$ of making a diameter is always meant with respect to the metric of the space $E$. Having a subset $A$ of a linear space $X$ we denote its convex hull (i.e., intersection of all convex sets of $X$ containing $A$) by $\co A$. Throughout the paper we also adopt a convention that $\inf \emptyset = \infty$.

A useful characterization of resolvable sets says that a set $H\subset X$ is resolvable if and only if there exist an ordinal $\kappa$, a subset $I$ of $[0,\kappa)$ and an increasing sequence of open sets $\emptyset=U_0\subset U_1\subset \cdots \subset U_\kappa=X$ satisfying
\begin{itemize}
	\item $U_\gamma = \bigcup\{U_\lambda:\lambda<\gamma\}$ for a limit ordinal $\gamma\leq\kappa$,\\
 	\item $H=\bigcup\{U_{\gamma+1}\setminus U_\gamma: \gamma\in I\}$.
\end{itemize}
A system of sets $\{U_{\alpha+1}\setminus U_\alpha: \alpha<\kappa\}$ with such properties is called a \emph{resolvable partition} of $X$.

Let $\F$ be a family of subsets of $X$. Then $f$ is \emph{$\ep$-$\sigma$-fragmented by sets of $\F$} if there is a system $\{X_n:n\in\en\}\subset\F$ covering $X$ such that $f|_{X_n}$ is $\ep$-fragmented for each $n\in\en$.

We set
\[
\sfrag_{\F} (f) = \inf \{\ep>0: f \text{ is $\ep$-$\sigma$-fragmented by sets of $\F$} \}
\]
and say that $f$ is \emph{$\sigma$-fragmented by sets of $\F$} if $\sfrag_{\F}(f) = 0$. For the sake of brevity we write just $\sfrag$ instead of $\sfrag_{\Hs(X)}$ if $\F = \Hs(X)$ (cf. \cite{jaynefrag}). The space of all $\sigma$-fragmented functions $f:X\to E$ (i.e. $\sfrag(f)=0$) will be denoted by $\sfra(X,\mathbb{R})$.

Let us remind that a mapping $f$ from a topological space $X$ to a metric space $E$ has the \emph{point of continuity property} (or briefly PCP) provided for every nonempty closed $F\subset X$, the restriction $f|F$ of $f$ to $F$ has a point of continuity.

Following Koumoullis in \cite{koumoullis}, we consider a certain generalization of functions of Baire class~1 (i.e. functions expressible as pointwise limits of sequences of continuous functions) from the metric spaces within the context of topological spaces. Let $\Bos(X)$ stand for the algebra generated by open sets and $\Hs(X)$ for the algebra of all resolvable sets. Let $\Sigma_2(\Bos(X))$ stand for the family of all countable unions of sets from $\Bos(X)$, analogously we define $\Sigma_2(\Hs(X))$. We say that $f:X\to E$ is a mapping of the \emph{first Borel class} if $f$ is $\Sigma_2(\Bos(X))$-measurable, that is, $f^{-1}(U)\in \Sigma_2(\Bos(X))$ for every open set $U\subset X$. The family of such mappings is designated by $\Bof_1(X,E)$. Similarly we define the family of all mappings of the \emph{first resolvable class} or \emph{first $H$-class} and denote them by $\Hf_1(X,E)$. If $X$ is metrizable then $\Hf_1(X,E)$ and $\Bof_1(X,E)$ coincide.

It is easy to check that any function of Baire class~1 is also of the first Borel class (see, e.g., \cite[Exercise 3.A.1]{fine}) and that any function of the first Borel class is also of the first $H$-class.


The aim of our paper is twofold. First, in Section~\ref{distances2mappings} of the paper we extend some results of \cite{ca-ma-ra}, \cite{ca-b1} and \cite{spurny-studia} by computing the distance of a general mapping to the family of mappings of the first resolvable class via the quantity $\frag$. In Section~\ref{quant-diff} we provide results analogous to those in \cite{quantitative}, \cite{ca-b1} concerning quantitative difference between countable compactness and compactness in $\Hf_1(X,E)$. Kindred results may be also found in \cite{angosto-dis}.

Second, in Section~\ref{osc-rank} we study a class of mappings with a countable oscillation rank and relate its basic properties to the aforementioned classes of mappings. This rank has been in a less general context considered by many authors (see, e.g., \cite{Haydon91oncertain}). It also possess a connection to the Szlenk index (e.g., \cite[Definition 8.5]{fabian2010banach}) in the following way. Assume $X$ is an infinite-dimensional Banach space and $K$ a $w^*$-compact subset of $X^*$. Let $f$ be the identity mapping of $(K,w^*)$ onto $(K,\|\cdot\|)$. Then the Szlenk index of $K$ is equal to the oscillation rank of $f$.

The investigation of the class of mappings with countable oscillation rank is also motivated by the following well known characterization of functions of Baire class~1: Let $E$ be a compact metrizable space. Then a function $f:E\to\er$ is of Baire class~1 if and only if its oscillation rank is countable (\cite[Proposition 2]{kechris1990classification}).


\section{Preliminaries}
This section contains auxiliary results dealing particularly with an approximation of fragmented mappings and provide us a necessary toolbox for Section~\ref{distances2mappings}.

\begin{lemma}\label{nastejno}
Let $X$ be a hereditarily Baire topological space, $E$ a metric space. Then for every mapping $f:X\to E$, $\sfrag(f) = \frag(f)$.
\end{lemma}

\begin{proof}
The inequality $\sfrag(f) \leq \frag(f)$ follows immediately from the definitions.

If $\sfrag(f) = \infty$ then the remaining inequality holds. Provided $\sfrag(f) \neq \infty$, let us assume that $\sfrag(f) < \ep$ for some $\ep\in\er$. Then there exists a sequence $\{H_n: n\in\en\} \subset \Hs(X)$ covering $X$ such that $f|_{H_n}$ is $\ep$-fragmented for each $n\in\en$. Given a closed subset $F$ of $X$ consider the resolvable sets $E_n=H_n\cap F$, $n\in\en$. 
From the proof of \cite[Proposition 2.1(iv)]{koumoullis} it follows that for each $n\in\en$ there exist sets $U_n$, $N_n$ respectively open and nowhere dense in $F$ satisfying $E_n=U_n\cup N_n$. Since $F$ is a Baire space there exists $k\in\en$ such that $U_k\neq\emptyset$.

Employing $\ep$-fragmentability of $f|_{H_k}$ we find an open subset $V$ of $X$ such that $V\cap U_k \neq \emptyset$ and $\diam f(V\cap U_k) < \ep$. Since $V\cap U_k$ is relatively open in $F$, $f$ is $\ep$-fragmented and $\frag(f) \leq \sfrag(f)$.

\end{proof}

\begin{lemma}\label{aprox}
Let $X$ be a topological space, $(E,\rho)$ a metric space and $f:X\to E$ an $\ep$-fragmented mapping for some $\ep>0$. Then there exists a mapping $g:X\to E$ which is constant on each set of some resolvable partition of $X$ and $\rho(f(x),g(x))<\ep$ for every $x\in X$.

Moreover, if $E=\er$ then such $g$ can be found so that $\rho(f(x),g(x))<\frac{\ep}2$ for every $x\in X$.
\end{lemma}

\begin{proof}
We find an ordinal $\Gamma$ and construct sets $G_\alpha$, $F_\alpha$ for $\alpha<\Gamma$ by transfinite induction. 

Set $G_0=F_0=\emptyset$. Let us assume that sets $G_\xi$, $F_\xi$ are constructed for every $\xi < \gamma$. If $\bigcup_{\xi<\gamma} G_\xi = X$ we set $\Gamma:=\gamma$ and stop the construction. Otherwise, due to $\ep$-fragmentability of $f$ there exists an open set $G_\gamma$ satisfying 
$$\diam f(G_\gamma \setminus \bigcup_{\xi<\gamma} G_\xi) < \ep.$$
We set $F_\gamma = G_\gamma \setminus \bigcup_{\xi<\gamma} G_\xi$.

Then $\{F_\gamma: \gamma<\Gamma\}$ is a resolvable partition of $X$ and we can define a mapping $g:X\to E$ as follows: Given $\gamma<\Gamma$ and $x\in F_\gamma$ we set $g(x) = t_\gamma$ where in the general case $t_\gamma\in f(F_\gamma)$ is chosen arbitrarily whereas in case $E=\er$ we take $t_\gamma$ as the center of $\co(f(F_\gamma))$.

Moreover, in the general case the inequality $\rho(f(x),g(x)) < \ep$ holds for every $x\in E$ and if $E=\er$ then even $\rho(f(x),g(x)) < \frac{\ep}2$ for every $x\in E$ holds.
\end{proof}

\begin{remark}\label{r:schod}
Realize that a mapping $g:X\to E$ from a topological space $X$ to a metric space $E$, which is constant on each set of some resolvable partition $\{G_{\alpha+1}\setminus G_\alpha : \alpha < \kappa \}$ of $X$, is fragmented. Indeed, let $\ep>0$ and $F$ be a closed subset of $X$ and set $\delta = \inf \{ \alpha \leq \kappa: G_\alpha\cap F \neq \emptyset \}$. It can be easily observed that $\delta$ is a successive ordinal, i.e., there exists $\gamma<\kappa$ such that $\delta=\gamma+1$. Hence, the mapping $g$ is constant on the set
$$F\cap G_{\gamma+1}= F\cap (G_{\gamma+1}\setminus G_\gamma)$$
and therefore $\diam g(F\cap G_{\gamma+1})=0$. The function $g$ is thus $\ep$-fragmented for every $\ep>0$ and hence fragmented.
\end{remark}

\begin{cor}
Let $X$ be a topological space and $(E,\rho)$ a metric space. Then $f:X\to E$ is a fragmented mapping if and only if there exists a sequence of mappings $(f_n)_n$ which converges uniformly to $f$ such that every $f_n$ is constant on each set of some resolvable partition of $X$.
\end{cor}

\begin{proof}
"Only if part" is a straightforward consequence of Lemma~\ref{aprox}.

For the proof of the "if part" let $\ep>0$. Then there exists $k\in\en$ such that $\rho(f(x),f_k(x))<\frac{\ep}3$ for every $x\in X$. Then, due to Remark~\ref{r:schod} for every closed set $F$ in $X$ there is an open set $U$ in $X$ satisfying $U\cap F\neq\emptyset$ and $\diam f_k(U\cap F)<\frac{\ep}3$. Hence, $\diam f(U\cap F) <\ep$ which concludes the proof. 

\end{proof}

\section{Relation of $\Hf_1(X,E)$ and $\sfra(X,E)$}


We start the section with a useful observation.

\begin{prop}
Let $X$ be a topological space and $(E,\rho)$ a metric space. Then $\sfra(X,E)\subset \Hf_1(X,E)$.
\end{prop}
\begin{proof}
Let $f\in\sfra(X,E)$. 

\end{proof}

We aim to show in Propositions~\ref{separable} and~\ref{suslinovaina} that, under an extra hypothesis, even an equality $\sfra(X,E)=\Hf_1(X,E)$ holds.

\begin{prop}\label{separable}
Let $X$ be a topological and $(E,\rho)$ a separable metric space. Then $\Hf_1(X,E)\subset\sfra(f,E)$.
\end{prop}

\begin{proof}
Given $\ep>0$, we find a set of points $\{y_i\in E: i\in\en\}$ such that the system of open balls $\{B(y_i,\frac{\ep}{2}): i\in\en\}$ covers $E$. Then $f^{-1}(B(y_i,\frac{\ep}{2}))\in \Sigma_2(\Hs(X))$, and thus there exists a family of resolvable sets $\{H_{i,n}$: $i,n\in\en\}$ satisfying

\[
f^{-1}(B(y_i,\frac{\ep}{2}))=\bigcup_{n=1}^\infty H_{i,n},\quad i\in \en.
\]

Then $\{H_{i,n}: i,n\in\en\}$ is a countable cover of $X$ such that $\diam f(H_{i,n})\leq\ep$, for every $i,n\in\en$. We may conclude that $\sfrag(f)=0$.
\end{proof}

To formulate the following proposition we need to recall some topological notions. A subset of a topological space $X$ is said to be \emph{Suslin} if it arises as a result of the Suslin operation applied to a system of closed sets in $X$. 

A topological space $X$ is an \emph{absolute Suslin space} if it is homeomorphic to a Suslin set in a complete metric space.

A family $\V$ of sets in a topological space $X$ is \emph{Suslin-additive} if $\bigcup \U$ is Suslin for every subfamily $\U \subset \V$. 

Let us also remind that a family $\V$ is called \emph{discrete} if every point $x\in X$ has a neighborhood that intersects at most one set of $\V$.
A family $\{A_\alpha:\alpha\in I\}$ is said to be \emph{$\sigma$-discretely decomposable} if there exists a family $\{A_{\alpha,n}: \alpha\in I, n\in\en\}$ of sets in $X$ such that $\{A_{\alpha,n}: \alpha\in I\}$ is discrete for every $n\in\en$ and $A_\alpha=\bigcup_{n\in\en} A_{\alpha,n}$.

For more details regarding these concepts consult, e.g., \cite{hansell1974characterizing}, \cite{hansell1998nonseparable}. Introducing the notion of an absolute Suslin space allows us to find another relation between mapping of the first $H$-class and the quantity $\sfrag$.

\begin{prop}\label{suslinovina}
Let $X$ be an absolute Suslin and $E$ a metric space. Then $\Hf_1(X,E)\subset\sfra(X,E)$.
\end{prop}

\begin{proof}

Let $\ep>0$. Stone's theorem (see \cite[Theorem 4.4.1]{engelking}) provides an open cover $\V$ of $E$ consisting of open sets of diameter smaller than $\ep$ such that $\V=\bigcup\V_n$, where a family $\V_n$ is discrete for every $n\in\en$.

Fix $n\in\en$. Then the family $\{ h^{-1}(V): V\in\V_n \}$ is disjoint and since any resolvable set in a metric space is Borel (see \cite[\S30, X, Theorem 5]{kuratowski-top}) it is also Suslin-additive.

Due to \cite[Theorem 2]{Ha}, the above-mentioned family is $\sigma$-discretely decomposable, namely, there exists a family $\{H_{V,k}: V\in \V_n, k\in\en\}$ such that the family $\{H_{V,k}: V\in \V_n\}$ is discrete for every $k\in\en$ and $h^{-1}(V)=\bigcup_{k\in\en} H_{V,k}$ for every $V\in\V_n$. Without loss of generality we may assume that $\{H_{V,k}: V\in \V_n, k\in\en\} \subset \Sigma_2(\Hs(X))$ (otherwise we would consider a family $\{ \overline{H_{V,k}}\cap h^{-1}(V): V\in\V_n\, k\in\en\}$ instead).

Then every $H_{V,k}$, where $V\in\V_n$, $k\in\en$, has a partition $\{H_{V,k,m}:m\in\en\}$ comprised of resolvable sets. We consider the sets
\[
H_{k, m, n} = \bigcup_{V\in \V_n} H_{V,k,m},\quad \text{where } k,m,n\in\en,
\]
which are unions of discrete families of resolvable sets and hence resolvable as well.

A straightforward verification suffices to realize that $X = \bigcup_{k,m,n} H_{k,m,n}$ and the mapping $h\r_{H_{m,n,k}}$ is $\ep$-fragmentable for each $k,m,n\in\en$. Hence $\sfrag(h) < \ep$ and $\sfrag(h) = 0$.

\end{proof}

Now, we provide examples showing that $\sfra(X,E)=\Hf_1(X,E)$ is not true in general.

\begin{example2}[{see \cite[Example 2.4(3)]{koumoullis}}]\label{e:kou}
Assuming Martin's Axiom and the negation of the continuum hypothesis, \cite[p. 162]{marsol} assures the existence of an uncountable set $Z$ in $\er$ whose every subset is an $F_\sigma$ set in $Z$. Let $X$ be the set $Z$ endowed with the Euclidean metric and $E$ be $Z$ with the discrete metric. Consider the identity mapping $h: X\to E$.

Since $X$ is a metric space, the family of sets $\Sigma_2(\Hs(X))$ corresponds to the family of $F_\sigma$ subsets of $X$ (see \cite[Proposition 3.4(d2)]{spurny-amh}). Hence $h\in\Hf_1(X,E)$.

On the other hand, $\sfrag(h)=1$. Indeed, given $0<\ep<1$ and a system $\{X_n:n\in\en\}$ of resolvable sets covering $X$, the mapping $h\r_{X_n}$ cannot be $\ep$-fragmented on each $X_n$. If that was the case, we would select $k\in\en$ with $X_k$ uncountable. Observe that every subset of $X_k$ would have an isolated point due to the $\ep$-fragmentability of $h\r_{X_k}$. Since $X_k$ has a countable base, $X_k$ would be a countable set, which would contradict the cardinality of the selected $X_k$. Since clearly $\sfrag(h)\leq 1$, the desired conclusion follows.
\end{example2}

\begin{remarks}
\begin{enumerate}[a)]
\item Let us remark that $X$ constructed in Example~\ref{e:kou} is not hereditarily Baire. Up to our knowledge, it is not known whether there exists a counterexample to $\sfra(X,E)=\Hf_1(X,E)$ in ZFC.
\item If we assume $X$ to be hereditarily Baire then $\sigma$-fragmentability is equivalent to fragmentability (see Lemma~\ref{nastejno}) and this is equivalent to having the PCP (see e.g. \cite[Proposition 2.]{koumoullis}). The equality $\sfra(X,E)=\Hf_1(X,E)$ is then equivalent the statement that a function of first H-class has the PCP.
\item By \cite[Theorem 4.1]{koumoullis}, if $X$ is a hereditarily $t$-Baire space (which is a subclass of hereditarily Baire spaces) and $E$ a metric space with cardinality less than the least ($\{0,1\}$-)measurable cardinal, then $f\in\Hf_1(X,E)$ if and only if $f$ has the PCP. \label{r:cons}

\item It is known that the axiom of constructibility ($V=L$) is independent on ZFC (see \cite{Shoenfield59}) and implies non-existence of measurable cardinals (see \cite{Scott61}). Hence (due to {\ref{r:cons})}), it is consistent with ZFC that $\sfra(X,E)=\Hf_1(X,E)$ for every hereditarily Baire space $X$ and metric space $E$.

\item As was shown by O. Kalenda in \cite{K97} and \cite{K98}, the cardinality restriction on the metric space $E$ in {\ref{r:cons})} is necessary. More precisely, if there is a ($\{0,1\}$-)measurable cardinal, then  there exist a (completely regular) hereditarily $t$-Baire space $X$ and a metric space $E$ such that $f\in\Hf_1(X,E)$ and $f$ has no point of continuity.

\item Also, by \cite[Examples 2.4(2)]{koumoullis}, assuming the existence of a real-valued measurable cardinal smaller or equal to $2^{\aleph_0}$, there exist a (Hausdorff) hereditarily Baire space $X$, a metric space $E$ and a mapping $h\in \Hf_1(X,E)$ such that $h$ has no point of continuity.
\end{enumerate}
\end{remarks}

\section{Distances to the function spaces $\sfra(X,E)$ and $\Hf_1(X,E)$}\label{distances2mappings}

The aim of this section is to estimate the distance between a general mapping from a topological space $X$ to a metric space $E$ and the families $\sfra(X,E)$ and $\Hf_1(X,E)$ via the quantity $\frag$. However, as is demonstrated later in this section, the results depend on the structure of spaces $X$ and $E$ in general.

Let $X$ be a topological space and $(E,\rho)$ a metric space. We define a distance of a pair of mappings $f,g:X\to E$ by 
\[
d(f,g)=\sup\{\rho(f(x),g(x)): x\in X\}.
\]
If $\F$ is a system of mappings from $X$ to $E$, we denote
\[
d(f,\F)=\inf\{d(f,g): g\in\F\}.
\]

\begin{prop}\label{dist-sfra}
If $X$ is a topological space, $(E,\rho)$ a metric space, then for every mapping $f:X\to E$ holds
\begin{equation}\label{dhaf}
\frac12\sfrag(f)\leq d(f, \sfra(X,E)) \leq \sfrag(f).
\end{equation}

Moreover, if $E=\er$, then
$d(f, \sfra(X,\er)) = \frac12 \sfrag(f)$.
\end{prop}

\begin{proof}
First, we prove the upper estimates $d(f, \sfra(X,E)) \leq \sfrag(f)$ and, for $f:X\to\er$, $d(f, \sfra(X,\er)) \leq \frac12 \sfrag(f)$.

If $\sfrag = \infty$ then it clearly holds. Now, suppose that $\sfrag(f)<\ep$ for some $\ep\in\er$. Then $X$ can be covered by a system $\{X_n: n\in\en\}$ of resolvable sets such that $f|_{X_n}$ is $\ep$-fragmented for every $n\in\en$.

If we construct a partition $\{Y_n:n\in\en\}$ of $X$ in a usual way 
\[
Y_1=X_1 \text{ and } Y_n = X_n \setminus \bigcup_{m=1}^{n-1} X_m,
\]
then $f|_{Y_n}$ is $\ep$-fragmented for every $n\in\en$.

For every $n\in\en$ we apply Lemma~\ref{aprox} and obtain a mapping $g_n: Y_n \to E$ constant on each set of some resolvable partition of $Y_n$ and satisfying $\rho(g_n(x),f(x)) <\ep$, $x\in Y_n$. Now we define a mapping $g:X\to E$ by setting $g(x) = g_n(x)$ for $x\in Y_n$, $n\in\en$.

Then $g$ is $\sigma$-fragmented and for every $x\in X$, it holds $\rho(f(x),g(x))<\ep$. Hence, $\rho(f,\sfra(X,E)) \leq \sfrag(f)$.

If $E=\er$, Lemma~\ref{aprox} provides a $\sigma$-fragmented approximation $g:X\to E$ such that $\rho(f(x),g(x))<\frac{\ep}2$, for every $x\in X$. Thus, $\rho(f,\sfra(X,E)) \leq \frac12 \sfrag(f)$.

\medskip

Second, we prove the lower estimate $\frac12 \sfrag(f) \leq d(f, \sfra(X,E))$.

If $d(f, \sfra(X,E)) = \infty$ then the statement is clearly valid. Otherwise, suppose $d(f, \sfra(X,E))<\alpha$ for some $\alpha\in\er$. Then there exists $h\in \sfra(X,E)$ such that $\rho(f(x),h(x)) < \alpha$ for every $x\in X$. Given an $\ep>0$, there exists a system $\{H_n:n\in\en\}$ of resolvable sets covering $X$ such that $h|_{H_n}$ is $\ep$-fragmented for each $n\in\en$.

Let $n\in\en$ and let $F$ be a closed subset of $H_n$. Then there exists an open subset $U$ of $X$ such that $U\cap F\neq \emptyset$ and $\diam h(U\cap F)<\ep$. From 
\[
\aligned
\diam f(U\cap F) &= \sup_{a,b\in U\cap F} \rho(f(a),f(b)) \\
 &\leq \sup_{a,b\in U\cap F} \rho(f(a),h(a)) + \rho(h(a),h(b)) + \rho(h(b),f(b)) \leq 2\alpha + \ep,
\endaligned
\]
we conclude that $\sfrag(f) \leq 2\alpha$ and this finally yields
\[
\frac12 \sfrag(f) \leq d(f, \sfra(X,E)).
\]
\end{proof}

We would like to use the quantity $\sfrag$ (or $\frag$) to get a lower estimate of $d(f,\Hf_1(X,E))$. We proceed in two steps. First, we show that a certain qualitative property of the space $\Hf_1(X,E)$ already guarantees a quantitative estimate. Second, we bring out examples of families of spaces $X$, $E$ with $\Hf_1(X,E)$ having the mentioned qualitative property.

%
%
%
%
%

The results from the previous two sections provide a following corollary (compare with \cite[Theorem 2.5 and Corollary 2.6]{ca-b1}).

\begin{thm}\label{souhrn}
Let $X$ be a topological space and $E$ a metric space. If one of the following conditions is satisfied
\begin{enumerate}
	\item[(i)] $X$ is an  absolutely Suslin space, and
	\item[(ii)] $E$ is a separable metric space,
\end{enumerate}
then for every $f: X\to E$ holds
\[
\frac12 \sfrag(f) \le d(f,\Hf_1(X,E))\le \sfrag(f).
\]
Moreover, if $E=\er$, then
\[
d(f,\Hf_1(X,\er)) = \frac12 \sfrag(f).
\]
\end{thm}

\begin{proof}
A consequence of Propositions~\ref{dist-sfra},~\ref{separable} and~\ref{suslinovina}.
\end{proof}

\section{Quantitative difference between compactness and countable compactness in $\Hf_1(X,E)$}\label{quant-diff}

We follow the line of reasoning which has appeared in \cite{quantitative} and \cite{ca-b1} where the quantitative differences between compactness and countable compactness in $C(X,E)$ (see \cite[Theorem 2.3]{quantitative}) and $B_1(X,E)$ (see \cite[Corollary 3.2]{ca-b1}) have been studied. The goal of this section is to provide analogous results for $\Hf_1(X,E)$.

For this purposes we adopt the following notions. Let $X$ be a topological space. Given a subset $A\subset X$, we denote the set of all sequences in $A$ by $A^\en$ and the set of all cluster points of a sequence $\varphi\in A^\en$ in $X$ by $\clust(\varphi)$.

Let $A$, $B$ be nonempty subsets of a metric space $(E,d)$. We employ a notion of a \emph{usual distance} between $A$ and $B$ defined by
\[
d(A,B) = \inf\{d(a,b): a\in A, b\in B\}
\]
and the \emph{Hausdorff non-symmetrical distance} from $A$ to $B$ defined by
\[
\hat{d} (A,B) = \sup\{d(a,B): a\in A\}.
\]

Let $S$ be a set and $X$ be a topological space. Then $(S^X,\tau_p)$ denotes a topological space generated by the pointwise convergence.

The following theorem is proved as \cite[Proposition 3.1]{ca-b1} provided $X$ is a separable metric space, nevertheless the same arguments serve equally well in a little more general setting:

\begin{prop}\label{clust}
Let $X$ be a second countable topological space, $E$ a metrizable space and $H$ a relatively compact subset of $(E^X,\tau_p)$. Then
\[
\sup_{f\in {\overline{H}}^{\tau_p}} \frag(f) = \sup_{\varphi\in H^\en} \inf\{ \frag(f): f\in \clust(\varphi) \}.
\]
\end{prop}

Let $X$ be a topological space, $(E,d)$ a metric space and $H$ a relatively compact subset of the space $(E^X,\tau_p)$. Then we may define the quantity
\[
\ck(H) = \sup_{\varphi\in H^\en} d(\clust(\varphi), \Hf_1(X,E)).
\]

\begin{remark}\label{clust-l}
Note that for a general topological space $X$, a metric space $E$ and a relatively compact set $H\subset (E^X,\tau_p)$, trivially
\[
\ck(H) \leq \hat{d}(\overline{H}^{\tau_p}, \Hf_1(X,E)).
\]

If $H$ is moreover a relatively countably compact subset of $(\Hf_1(X,E),\tau_p)$, then $\ck(H) = 0$ (use, e.g., \cite[Theorem 3.10.3]{engelking}).
\end{remark}

Now, with the aid of the previous results we deduce a few corollaries.

\begin{cor}\label{clust-2}
Let $X$ be a second countable hereditarily Baire topological space, $E$ be a metric space and $H$ be $\tau_p$-relatively compact subset of $E^X$. If moreover one of the following conditions is satisfied
\begin{enumerate}[(i)]
\item $X$ is an absolute Suslin topological space,
\item $E$ is a separable metric space,
\end{enumerate}
then
\[
\ck(H) \leq \hat{d}(\overline{H}^{\tau_p}, \Hf_1(X,E)) \leq 2\ck(H).
\]
\end{cor}
\begin{proof}
The inequality $\ck(H) \leq \hat{d}(\overline{H}^{\tau_p}, \Hf_1(X,E))$ was explained in Remark~\ref{clust-1} already.

Now, we prove the inequality $\hat{d}(\overline{H}^{\tau_p}, \Hf_1(X,E)) \leq 2\ck(H)$. Taking the definition of $\hat{d}$ into consideration and following Theorem~\ref{souhrn} together with Lemma~\ref{nastejno} we infer that
\[
\aligned
\hat{d}(\overline{H}^{\tau_p}, \Hf_1(X,E)) = \sup_{f\in \overline{H}^{\tau_p}}d(f, \Hf_1(X,E)) \leq \sup_{f\in \overline{H}^{\tau_p}} \sfrag(f) = \sup_{f\in \overline{H}^{\tau_p}} \frag(f) = (\ast).
\endaligned
\]
We continue by employing Proposition~\ref{clust} and the closing argument is again due to Lemma~\ref{nastejno} and Theorem~\ref{souhrn}
\[
\aligned
(\ast) &= \sup_{\varphi\in H^\en} \inf\{ \frag(f): f\in \clust(\varphi) \} \\
		&= \sup_{\varphi\in H^\en} \inf\{ \sfrag(f): f\in \clust(\varphi) \} \\
			&\leq \sup_{\varphi\in H^\en} \inf\{2 d(f, \Hf_1(X,E)): f\in \clust(\varphi) \} = 2\ck(H).
\endaligned
\]
\end{proof}

\begin{cor}
Let $X$ be a separable completely metrizable space and $H$ be a $\tau_p$-relatively compact subset of $\er^X$. Then
\[
\hat{d}(\overline{H}^{\tau_p}, \Hf_1(X,\mathbb{R})) =\ck(H).
\]
\end{cor}
\begin{proof}
The argument goes along similar pattern as the proof of Corollary~\ref{clust-2}, using Theorem~\ref{souhrn} for the case $E=\mathbb{R}$.
\end{proof}

\begin{remark}
Let $X$ be a topological space, $E$ a metric space and $H$ a relatively compact subset of $(E^X,\tau_p)$. Note that:
\begin{itemize}
	\item If $E$ is a relatively countably compact subset of $H_1(X,Z),\tau_p)$, then $\ck(H)=0$.
	\item If $\hat{d}(\overline{H}^{\tau_p}, \Hf_1(X,E))=0$, then $H$ is a relatively compact subset of $\Hf_1(X,E)$.	
\end{itemize}
Assuming that the hypothesis of Corollary~\ref{clust-2} is satisfied, it follows that in $(H_1(X,Z),\tau_p)$ the relatively countably compact sets are also relatively compact. This is (after results in \cite{ca-b1}) another generalisation of the classical Rosenthal's result  (\cite{Ro74}) stating that in $(B_1(X,\mathbb{R}),\tau_p)$, where $X$ is a Polish space, the relatively countably compact sets are relatively compact.
\end{remark}

\section{Oscillation rank}\label{osc-rank}
In this section we recall a definition of the oscillation rank which has been adopted by many authors, including S.~Argyros, R.~Haydon, A.~S.~Kechris and A.~Louveau (see, e.g., \cite{Haydon91oncertain}, \cite{kechris1990classification}). However, this rank has been so far investigated for functions defined on metrizable compact spaces. The main purpose of this section is to provide a view of the situation considering the oscillation rank for mappings from topological spaces to metric spaces.

We adhere to a standard convention that $\inf \emptyset = \infty$ and $\infty$ is greater then any ordinal.
\begin{definition}
Let $X$ be a topological space, $(E,\rho)$ a metric space and $\ep>0$. For a given mapping $f:X\to E$, we construct for every ordinal $\alpha$ an open se $U_\alpha$. We proceed by transfinite induction. 

Let $U_0 = \emptyset$. Assume that $U_\gamma$ is constructed for every ordinal $\gamma<\alpha$. If $\alpha = \gamma'+1$ is a successive ordinal then we set 
\[
U_{\alpha} = U_{\gamma'} \cup \{ x\in X\setminus U_{\gamma'}: \text{$\exists$ open set $U \ni x$ such that } \diam f(U\setminus U_{\gamma'}) < \ep \}.
\]
If $\alpha$ is a limit ordinal we set $U_\alpha = \bigcup_{\gamma<\alpha} U_\gamma.$

We define $\orank(f,\ep)$ as the first ordinal $\alpha$ satisfying $U_\alpha = X$ and if such an ordinal does not exist we set $\orank(f,\ep) = \infty$.

Then we define the \emph{oscillation rank of a mapping $f$} as 
\[
\orank (f) = \sup_{\ep>0} \orank(f,\ep).
\]

Further, we define $\orank^{*}(f,\ep)$ as the first ordinal $\kappa$ such that there exists a transfinite sequence $(V_\alpha)_{\alpha\leq\kappa}$ of sets in $X$ with the property $*(f,\ep)$, where we say that a transfinite sequence $(V_\alpha)_{\alpha\leq\kappa}$ of sets in $X$ has the \emph{property $*(f,\ep)$} if
\begin{itemize}
\item[(i)] The transfinite sequence $(V_\alpha)_{\alpha\leq\kappa}$ is nondecreasing, composed of open sets in $X$ and such that $V_0=\emptyset$ and $V_\kappa = X$.
\item[(ii)] For every limit ordinal $\gamma\leq \kappa$, $V_\gamma = \bigcup_{\alpha<\gamma} V_\alpha$.
\item[(iii)] For every ordinal $\alpha<\kappa$ and $x\in V_{\alpha+1}\setminus V_\alpha$,there exists an open neighborhood $V$ of $x$ such that $\diam f(V\setminus V_\alpha) < \ep$.
\end{itemize} 
If there is no such a transfinite sequence we set $\orank^{*}(f,\ep) = \infty$.

We define
\[
\orank^{*} (f) = \sup_{\ep>0} \orank^{*}(f,\ep).
\]

\end{definition}

\begin{remark}\label{uzavrene}
In some cases it is more convenient to use a slightly modified definition of $\orank(f,\ep)$ where instead of constructing open subsets $U_\alpha\subset X$ for ordinals $\alpha$ we construct their complements $F_\alpha$.
Then $F_0=X$, for a limit ordinal $\alpha$ we set $F_\alpha = \bigcap_{\beta<\alpha} F_\beta$, and for a successive ordinal $\alpha = \beta'+1$ we set

\[
\aligned
F_{\alpha} &= \{ x\in F_{\beta'}: \text{$\forall$ open set $U \ni x$, } \diam f(U \cap F_{\beta'}) \geq \ep \} \\
					 &= F_{\beta'}\setminus \bigcup \{U\subset X: \text{$U$ is open, }\diam f(U\cap F_{\beta'}) < \ep\}.
\endaligned					
\]

Then $\orank(f,\ep)$ is the first ordinal $\alpha$ such that $F_\alpha=\emptyset$. If such an ordinal does not exist we set $\orank(f,\ep)=\infty$.
\end{remark}

\begin{remark}
We have adopted the notation of \cite{Haydon91oncertain}, though our definition of the oscillation rank is formally different. A glimpse of these two definitions shows, however, that they provide identical concepts.
\end{remark}

The following lemma proves that concepts $\orank$ and $\orank^{*}$ actually coincide.

\begin{lemma}\label{nastejno2}
Let $X$ be a topological space, $(E,\rho)$ a metric space and $\ep>0$. Then for every mapping $f:X\to E$, $\orank(f,\ep) = \orank^{*}(f,\ep)$.
\end{lemma}

\begin{proof}
The inequality $\orank^{*}(f,\ep)\leq\orank(f,\ep)$ follows immediately from the definitions.

The reverse inequality is obvious if $\orank(f,\ep)=\infty$. So let $\orank(f,\ep)=\kappa\neq\infty$ and let $(V_\alpha)_{\alpha\leq\kappa}$ be a transfinite sequence of sets in $X$ with $*(f,\ep)$. We find a transfinite sequence $(U_\alpha)_{\alpha\leq\kappa}$ provided by the construction in the definition of $\orank(f,\ep)$. Now we observe that for every $\alpha\leq \kappa$, $V_\alpha \subset U_\alpha$.

Indeed, for $\alpha=0$ the statement follows immediately from the definitions. Assume that $\gamma\leq\kappa$ is an ordinal and the statement is valid for every $\alpha<\gamma$. If $\gamma$ is a limit ordinal then $V_\gamma = \bigcup_{\alpha<\gamma} V_\alpha \subset \bigcup_{\alpha<\gamma} U_\alpha = U_\gamma$. If $\gamma$ is a successive ordinal then $\gamma=\alpha+1$ for some ordinal $\alpha$ and
\[
U_{\alpha+1} = U_\alpha \cup \{ x\in X\setminus U_\alpha: \text{$\exists$ open set $U \ni x$ such that } \diam f(U\setminus U_\alpha) < \ep \}.
\]

Let $x\in V_{\alpha+1}$. Then either $x\in U_{\alpha}\subset U_{\alpha+1}$, or $x\in V_{\alpha+1}\setminus U_{\alpha}\subset V_{\alpha+1}\setminus V_\alpha$ and due to $*(f,\ep)$ there exists an open neighborhood $V$ of $x$ such that $\diam f(V\setminus V_\alpha) <\ep$ and thus also $x\in U_{\alpha+1}$.

It follows that $V_\kappa \subset U_\kappa$, where $V_\kappa=X$, and hence $\orank(f,\ep)\leq\kappa=\orank^*(f,\ep)$ which concludes the proof.
\end{proof}
 
We may ask whether the transfinite process described in the definition of the oscillation rank must always stop. A quite natural answer is given by the following lemma.

\begin{lemma}\label{nutnost}
Let $X$ be a topological space, $(E,\rho)$ a metric space and $f:X\to E$. Then $\orank (f) \neq \infty$ if and only if $f$ is fragmented.
\end{lemma}

\begin{proof}
If $f$ is fragmented then a moment of reflection shows that $\orank (f) \neq \infty$. Specifically, $\orank(f)\leq\card(X)$, where $\card(X)$ is a cardinality of $X$.

Suppose $f$ is not fragmented. Then there exists a nonempty closed set $F\subset X$ and $\ep>0$ such that for every open set $U\subset X$ intersecting $F$, $\diam f(F\cap U) > \ep$. We prove that then $\orank (f,\ep) = \infty$. Suppose the contrary, let $\orank(f,\ep)=\kappa\neq\infty$ and let $(U_\alpha)_{\alpha\leq\kappa}$ be the transfinite sequence constructed in the definition of $\orank(f,\ep)$.

Let $\delta = \min \{\alpha: U_\alpha \cap F\neq\emptyset\}$. It is easy to see that $\delta\leq\kappa$ is a successive ordinal, hence there exists an ordinal $\gamma$ such that $\delta=\gamma+1$.

We choose a point $x\in F\cap U_{\gamma+1} = (F\cap U_{\gamma+1})\setminus U_\gamma$. Therefore, we can find an open neighborhood $U$ of $x$ such that $\diam f(U\setminus U_\gamma) < \ep$. Then
\[
\emptyset\neq U\cap F \cap U_{\gamma+1} = U\cap F\cap (U_{\gamma+1}\setminus U_\gamma) \subset U\setminus U_\gamma
\]
and so $\diam f(U\cap F \cap U_{\gamma+1})<\ep$ which contradicts the properties of the set $F$. We conclude that $\orank(f,\ep) = \infty$.
\end{proof}

The following lemma shows that the oscillation rank $\orank$ may attain a value of arbitrary ordinal less or equal $\omega_1$.

\begin{lemma}\label{prostory}
For every ordinal $\alpha$ there exists a completely metrizable topological space $X_\alpha$ and a function $f_\alpha:X_\alpha\to \{0,1\}$ satisfying $\orank(f_\alpha) = \alpha$ and $f_\alpha\in\Bof_1(X,E)$.
\end{lemma}

\begin{proof}
Let $\alpha$ be an ordinal. Then we find a completely metrizable topological space $X_\alpha$ and its subset $A_\alpha\subset X_\alpha$ such that for its characteristic function $\chi_{A_\alpha}$ holds $\orank(\chi_{A_\alpha}) = \alpha$.

The construction proceeds by transfinite induction on ordinal $\alpha$. Let $X_0 = A_0 =\emptyset$. Then for every $0<\ep<1$ holds $\orank(\chi_{A_0}) = \orank(\chi_{A_0},\ep) = 0$.

Assume that $\alpha<\omega_1$ is a limit ordinal and for every $\beta<\alpha$ we have constructed a completely metrizable topological space $X_\beta$ and its subset $A_\beta \subset X_\beta$ such that $\orank(\chi_{A_\beta}) = \beta$. We define topological spaces

\[
X = \bigcup_{\beta<\alpha} X_{\beta}\times\{\beta\} \enskip \mbox{ and } \enskip A = \bigcup_{\beta<\alpha} A_{\beta}\times\{\beta\}
\]

equipped with the disjoint union topology (e.g., \cite[Section~2.2]{engelking}). A space $X$ is completely metrizable due to \cite[p.~270]{engelking}. 

Now, let us fix $0<\ep<1$. For $\beta<\alpha$ let $(U^{\gamma}_\beta)_\gamma$ be the open sets emerging in the definition of $\orank(\chi_{A_\beta}, \ep)$. The system of sets
\[
\big\{ \bigcup_{\beta<\alpha } U^{\gamma}_\beta \times \{\beta\}: \gamma\leq\orank(\chi_A,\ep)\big\}
\]
is then precisely that one constructed in the definition of $\orank(\chi_A,\ep)$ and $\gamma=\alpha$ is clearly the first ordinal satisfying $\bigcup_{\beta<\alpha } U^{\gamma}_\beta \times \{\beta\} = X$. Hence, $\orank(\chi_A) = \alpha$ and we may set $X_\alpha=X$ and $A_\alpha=A$.

Let us assume that for an ordinal $\alpha$ we have a corresponding set $A_\alpha$ and a topological space $X_\alpha$ constructed. We define topological spaces
\[
X= \{p\} \cup \big(\bigcup_{n\in\en} X_\alpha \times \{n\}\big)  \enskip \mbox{ and } \enskip A= \{p\} \cup \big(\bigcup_{n\in\en} A_\alpha \times \{n\}\big)
\]
with the following topology: points of the form $(x,n)\in X_\alpha\times \{n\}$, where $n\in\en$, has a basis of neighborhoods consisting of sets $U\times \{n\}$ where $U\subset X_\alpha$ is a neighborhood of $x\in X_\alpha$ and the point $p\in X$ has a basis of neighborhoods consisting of sets of the form
\[
\{p\} \cup \bigcup_{n\geq n_0} X_\alpha \times \{n\},
\]
where $n_0\in\en$. The described topological space $X$ is completely metrizable. Indeed, let $\rho_n$ be a compatible complete metric $X_\alpha$ bounded by $\frac1n$ (see \cite[Proposition~4.3.8]{engelking}). Then we can define compatible complete metric $\rho$ on $X$ by
\[
\rho(x\times\{n\}, y\times\{m\}) =
\begin{cases}
\rho_n(x,y), & \text{ if } m=n\in\en, x,y\in X_\alpha,\\
\max\{\frac1n,\frac1m\}, & \text{ if } m\neq n, x,y\in X_\alpha
\end{cases}
\]
and
\[
\rho(x\times\{n\},p) = \frac1n.
\]
It can be easily verified that $\rho$ is a compatible metric. We will only prove the triangular inequality $\rho(x,z)+\rho(z,y)\geq \rho(x,y)$: 
\[
x,y\in X\times\{n\} \text{, }
\begin{cases}
z\in X_\alpha\times\{n\}:& \text{ Clear, since $(X_\alpha,\rho_n)$ is complete metric space.}\\
z\in X_\alpha\times\{m\}, m\neq n:& \rho(x,z)+\rho(z,y) =2\max\{\frac1n,\frac1m\} \geq \frac1n \geq \rho(x,y),\\
z=p:& \rho(x,z)+\rho(z,y) = \frac1n + \frac1n \geq \frac1n =\rho(x,y).
\end{cases}
\]
If $x\in X\times\{n\}$, $y\in X\times\{m\}$, $n,m\in\en$, $m>n$, then 
\[
\rho(x,z)+\rho(z,y)
\begin{cases}
 \geq \frac1n=\rho(x,y),& \text{ if } z\in X_\alpha\setminus\{p\},\\
  = \frac1n+\frac1m \geq \rho(x,y),& \text{ if } z=p.
\end{cases}
\]
If $x=p$ and $y\in X_\alpha\times\{n\}$, $z\in X_\alpha\times\{m\}$ then 
\[
\rho(x,z)+\rho(z,y) = \frac1m + \max{\frac1n,\frac1m} \geq \frac1n=\rho(x,y).
\]

Completeness of the metric $\rho$ is also easy to realize. Having a Cauchy sequence $(x_n)\subset X$ then either there exists $n_0,m\in\en$ such that $(x_n)_{n\geq n_0}\subset X_\alpha\times\{m\}$ (then $(x_n)$ converges to an element of $X_\alpha\times\{m\}$) or for every $m\in\en$ there exists $n_0\in\en$ such that $(x_n)_{n\geq n_0}\subset \{p\} \cup \big(\bigcup_{n\geq m} X_\alpha \times \{n\}\big)$ (then $(x_n)$ converges to $p$).

Following the inductive process in the definition of $\orank(f,\ep)$ (considering Remark~5.5.2) we construct closed sets $\{F_\beta: \beta \leq \orank(f,\ep)\}$. 

Our aim is to prove that $F_{\alpha} = \{p\}$. Once this assertion is justified it follows that $\orank(\chi_A)= \orank(\chi_A,\ep)= \alpha+1$ and we may set $X_{\alpha+1} = X$ and $A_{\alpha+1} = A$.

First, we take a system $\{H_\beta: \beta\leq\alpha\}$ of closed sets $\{H_\beta: \beta\leq\alpha\}$ from the definition of $\orank(\chi_{A_\alpha},\ep)$ (considering Remark~5.2.2). It is clear that $F_\beta \cap (X_\beta\times\{n\}) = H_\beta\times\{n\}$, $\beta\leq\alpha$, $n\in\en$. Since $F_\beta$, $\beta\leq\alpha$, is a closed set we can infer that $F_\beta = \{p\} \cup \bigcup_{n\in\en} H_\beta\times\{n\}$ for every $\beta\leq\alpha$. Specially,
\[
F_\alpha = \{p\} \cup \bigcup_{n\in\en} H_\alpha\times\{n\} = \{p\},
\]
which was to prove.

\end{proof}

\begin{remark}
According to \cite[Proposition 2.8]{Haydon91oncertain} for all $\alpha<\omega_1$ there exists a quasireflexive (of order 1) Banach space $Q_\alpha$ such that $Q_\alpha^{**}=Q_\alpha \oplus \langle f_\alpha\rangle$ where $\beta(f_\alpha)>\alpha$.
\end{remark}

Let $X$ be a topological space and $E$ a metric space. We denote the set of the mappings $f:X\to E$ satisfying $\orank(f)<\omega_1$ by the symbol $\oorank(X,E)$.

Next we prove the stability of $\oorank(X,E)$ under taking uniform limits of nets (Lemma~\ref{ulimity}) and composing with a uniformly continuous mappings (Lemma~\ref{uniform}).

\begin{lemma}\label{ulimity}
Let $X$ be a topological space and $(E,\rho)$ a metric space. If $(f_\gamma)_{\gamma\in\Gamma}$ is a net of mappings in $\oorank(X,E)$ which converges uniformly to $f$, then $f\in\oorank(X,E)$.
\end{lemma}

\begin{proof}
Let $\ep>0$. Then there exists $\gamma\in\Gamma$ such that $\rho(f(x),f_\gamma(x))\leq\frac{\ep}3$ for every $x\in X$. We set $\kappa=\orank^{*}(f,\frac{\ep}3)<\omega_1$ and find a family $(U_\alpha)_{\alpha\leq\kappa}$ with $*(f_\gamma,\frac{\ep}3)$. Then for each $\alpha<\kappa$ and $x\in U_{\alpha+1}\setminus U_\alpha$ there exists an open neighborhood $U$ of $x$ such that $\diam f_\gamma(U\setminus U_\alpha) \leq\frac{\ep}3$ and so
\[
\diam f(U\setminus U_\alpha) \leq \diam f_\gamma(U\setminus U_\alpha) + \frac{2\ep}3 < \ep.
\]
Hence, $\orank(f,\ep) = \orank^{*}(f,\ep) \leq \kappa$ (see Lemma~\ref{nastejno2}), and after passing to the supremum over $\ep>0$ also $\orank(f) \leq \kappa <\omega_1$ and finally $f\in \oorank(X,E)$.
\end{proof}

\begin{lemma}\label{uniform}
Let $X$ be a topological space, $(E,\rho)$ and $(F,\sigma)$ metric spaces, $h:E\to F$ a uniformly continuous mapping. Then for every $f:X\to E$, $\orank(h\circ f) \leq \orank(f)$ (and hence if $f\in\oorank(X,E)$ then $h\circ f\in\oorank(X,F)$).

\end{lemma}

\begin{proof}
Let $\ep>0$. We find $\delta>0$ such that for any $x,y\in E$ satisfying $\rho(x,y)<\delta$, $\sigma(h(x),h(y))<\ep$. Lemma~\ref{nastejno2} allows us to find a family $\{U_\alpha:\alpha\leq\orank(f)\}$ with $*(f,\delta)$.

For each $\alpha<\orank(f)$ and $x\in U_{\alpha+1}\setminus U_\alpha$ there exists an open neighborhood $U$ of $x$ such that $\diam f(U\setminus U_\alpha)<\delta$. Hence
\[
\diam (h\circ f)(U\setminus U_\alpha) \leq \ep,
\]
which finalizes the proof of the statement.
\end{proof}

The subsequent example indicates that continuous mappings do not generally preserve the oscillation rank.
\begin{example}
Let $f:\qe^{+}\to \er^{+}$ be defined as $f(\frac{p}{q}) = \frac1q$ for $\frac{p}{q}\in\qe^{+}$ where $p,q\in\en$ are coprime numbers. Further, we define a function $h:\er^{+}\to\er^{+}$ by setting $h(x) = \frac1{x}$ for $x\in\er^{+}$. Then $h$ is a continuous function and $\orank(f)=2<\omega_1$. Yet, the composition $h\circ f:\qe^{+}\to\er^{+}$ is not fragmented and hence $\orank(h\circ g) = \infty$. 
\end{example}

\begin{proof}
Let $0<\ep<1$ a $M:=\{ \frac{p}{q}:p,q\in\en, q<q_0\}$. We calculate $\oorank(f)$ from the definition. Clearly, $M\neq\emptyset$ and $U_0=\emptyset$. Then $U_1=\qe^+\setminus M\neq\qe^+$ and $U_2=\qe^+$, since $M$ is discrete in $\qe^+$. Hence, $\orank(f)=2$.

The function $h\circ f$ is not fragmented because it is unbounded on every nonempty open subset of $\qe^+$.

\end{proof}

Lemmas~\ref{vector},~\ref{lattice} and~\ref{algebra} provide us some information about the stability of $\oorank(X,E)$ under taking various algebraic operations. First, we prove the following auxiliary result.

\begin{lemma}\label{diagonal}
Let $X$ be a topological space, $(E,\rho)$ a metric space, $\kappa_1,\kappa_2$ ordinals and $f,g:X\to E$ mappings satisfying $\orank(f)\leq\kappa_1$, $\orank(g)\leq \kappa_2$, where $\kappa_1,\kappa_2>0$. Then for diagonal mapping $f\Delta g:X\to (E\times E, \rho_{\max})$ defined by $(f\Delta g)(x) = (f(x),g(x))$, $x\in X$ satisfies $\orank(f\Delta g) \leq \kappa_1 \kappa_2$ ($\rho_{\max}$ is the maximal metric).
\end{lemma}
\begin{proof}
Let $\ep>0$. Due to Lemma~\ref{nastejno2} we may find $(U_\alpha)_{\alpha\leq\kappa_1}$ ibn $X$ with $*(f,\ep)$ and $(V_\alpha)_{\alpha\leq\kappa_2}$ in $X$ with $*(g,\ep)$.

We remind that considering $\gamma<\kappa_1\kappa_2$ there exists a uniquely determined pair of ordinals $\alpha,\beta$ such that $\gamma = \kappa_1\beta + \alpha$ and $\alpha<\kappa_1$, $\beta<\kappa_2$ (e.g., \cite[Theorem VII.5.4]{kuratowski-set}). Hence, we may define
\[
\aligned
G_\gamma = (U_\alpha\cap V_{\beta+1})\cup V_\beta, \gamma<\kappa_1\kappa_2, \text{ and } G_{\kappa_1\kappa_2} = X.
\endaligned
\]

In the remainder of the proof we verify that the transfinite sequence $(G_{\gamma})_{\gamma\leq\kappa_1\kappa_2}$ satisfies $*(f\Delta g,\ep)$. 
First, $G_\gamma$ is an open set in $X$ for every $\gamma\leq\kappa_1\kappa_2$, $G_0=\emptyset$ and $G_{\kappa_1\kappa_2} = X$. Let $\gamma,\gamma'<\kappa_1\kappa_2$ be ordinals satisfying $\gamma<\gamma'$. Then $\gamma=\kappa_1\beta + \alpha$ and $\gamma'=\kappa_1\beta' + \alpha'$ where either $\beta<\beta'$, or $\beta=\beta'$ together with $\alpha<\alpha'$. If $\beta<\beta'$, then
\[
G_\gamma = (U_\alpha\cap V_{\beta+1})\cup V_\beta \subset V_{\beta'} \subset G_{\gamma'}
\]

and if $\beta=\beta'$ with $\alpha<\alpha'$ is the case, then
\[
G_\gamma = (U_\alpha\cap V_{\beta+1})\cup V_\beta \subset (U_{\alpha'}\cap V_{\beta+1})\cup V_\beta = G_{\gamma'}.
\]

Second, let $\gamma\leq\kappa_1\kappa_2$ be a limit ordinal with the unique decomposition $\gamma=\kappa_1\beta + \alpha$, $\alpha<\kappa_1$, $\beta\leq\kappa_2$. Then either $\alpha>0$ and $\alpha$ is a limit ordinal:
\[
G_\gamma = (U_\alpha\cap V_{\beta+1})\cup V_\beta = \bigcup_{\alpha'<\alpha} \big((U_{\alpha'}\cap V_{\beta+1})\cup V_\beta\big) = \bigcup_{\gamma'<\gamma} G_{\gamma'},
\]
or $\alpha=0$ and $\kappa_1$ is a limit ordinal:
\[
G_\gamma = V_\beta = \bigcup_{\beta'<\beta} (V_{\beta'+1}\cup V_{\beta'}) = \bigcup_{\beta'<\beta} \bigcup_{\alpha'<\kappa_1}(U_{\alpha'}\cap V_{\beta'+1})\cup V_{\beta'} = \bigcup_{\gamma'<\gamma} G_{\gamma'},
\]
or $\alpha=0$ and $\beta$ is a limit ordinal:
\[
G_\gamma = V_\beta = \bigcup_{\beta'<\beta} V_{\beta'} = \bigcup_{\alpha'<\kappa_1} \bigcup_{\beta'<\beta} (U_{\alpha'}\cap V_{\beta'+1})\cup V_{\beta'} = \bigcup_{\gamma'<\gamma} G_{\gamma'}.
\]

Third, let $x\in G_{\gamma+1}\setminus G_{\gamma}$, where $\gamma = \kappa_1\beta + \alpha<\kappa_1\kappa_2$. A decomposition of $\gamma+1$ is then either $\kappa_1\beta + (\alpha+1)$ or $\kappa_1(\beta+1)$ (when $\alpha+1 = \kappa_1$). Both cases imply
\[
G_{\gamma+1}\setminus G_{\gamma} = (U_{\alpha+1}\setminus U_\alpha)\cap (V_{\beta+1}\setminus V_\beta).
\]
Hence if $x\in G_{\gamma+1}\setminus G_{\gamma}$, we can find open sets $U, V\subset G_{\gamma+1}$ enjoying the following properties: $x\in U$, $x\in V$, $\diam f(U\setminus U_\alpha)<\ep$ and $\diam g(V\setminus V_\beta)<\ep$. It follows that $\diam (f\Delta g)((U\cap V) \setminus G_\gamma)<\ep$.

Hence, $\orank^{*} (f\Delta g,\ep) \leq \kappa_1\kappa_2$ and employing Lemma~\ref{nastejno2} the proof is complete.

\end{proof}

\begin{lemma}\label{vector}
Let $X$ be a topological space and $E$ a Banach space. If $f,g\in\oorank(X,E)$ then $f+g\in\oorank(X,E)$.
\end{lemma}

\begin{proof}
Let $f,g\in\oorank(X,E)$. We define $h:E\times E \to E$ as $h(x,y)=x+y$ for every $x,y\in E$. It can be verified straightforwardly that the mapping $h$ is uniformly continuous and that $f+g = h \circ (f\Delta g)$. Applying Lemmas~\ref{diagonal} and~\ref{uniform} we get the required inequality $\orank(f+g)<\omega_1$.
\end{proof}

\begin{lemma}\label{lattice}
Let $X$ be a topological space and $E$ a Banach lattice. If $f,g\in\oorank(X,E)$ then $\inf(f,g),\sup(f,g)\in\oorank(X,E)$.
\end{lemma}

\begin{proof}
First, we realize that the lattice mappings $i:E\times E\to E$ and $j: E\times E\to E$, defined as $i(x,y) = \inf(x,y)$ and $j(x,y)=\sup(x,y)$ for $x,y\in E$, are uniformly continuous (cf. \cite[Chapter II, Proposition 5.2]{lattices}). For any $f,g\in \oorank(X,E)$ and $x\in X$ the following identities clearly hold:
\[
\aligned
\sup(f,g)(x,y)=\sup(f(x),g(x)) = (j\circ (f\Delta g))(x),\\
\inf(f,g)(x,y) =\inf(f(x),g(x)) = (i\circ (f\Delta g))(x).
\endaligned
\]
Thanks to Lemmas~\ref{diagonal} and \ref{uniform}, $\sup(f,g),\inf(f,g)\in\oorank(X,E)$.
\end{proof}

Concerning stability under taking products we were able to achieve only a partial result.
\begin{lemma}\label{algebra}
Let $X$ be a topological space, $E$ a commutative Banach algebra and $f,g\in\oorank(X,E)$ bounded mappings. Then $fg\in\oorank(X,E)$.
\end{lemma}
\begin{proof}
Since $fg = \frac12\left((f+g)^2 - f^2 - g^2\right)$, by Lemma~\ref{vector} it suffices to prove that for every bounded $f:X\to E$ such that $f\in\oorank(X,E)$, $f^2\in\oorank(X,E)$.

We may write $f^2$ as a composition $f^2 = \varphi \circ f$ where $\varphi: E\to E$ is defined as $\varphi(x)=x^2$ for each $x\in E$. Since $f$ is a bounded mapping, $\varphi|_{\range f}$ is uniformly continuous and employing Lemma~\ref{uniform} the proof is complete.
\end{proof}

However, the general problem remains unsolved.
\begin{ques} 
Let $X$ be a topological space, $E$ a commutative Banach algebra and $f,g\in\oorank(X,E)$. Is it then true that $fg\in\oorank(X,E)$?
\end{ques}

We summarize the previous stability results in the following proposition.

\begin{prop}
Let $X$ be a topological space. 
	\begin{itemize}
		\item[(a)] If $E$ is a Banach space then the space of all bounded elements of $\oorank(X,E)$ endowed with the supremum norm (i.e., $\|f\|=\sup\{\|f(x)\|:x\in X\}$) is a Banach space.
		\item[(b)] If $E$ is a Banach lattice then the space of all bounded elements of $\oorank(X,E)$ is a Banach lattice in a supremum norm and the pointwise ordering.
		\item[(c)] If $E$ is a commutative Banach algebra then the space of all bounded elements of $\oorank(X,E)$ a commutative Banach algebra.
	\end{itemize}
\end{prop}

\begin{proof}
A straightforward verification using Lemmas~\ref{ulimity},~\ref{vector},~\ref{lattice} and~\ref{algebra}.
\end{proof}

Now we clarify the relations between $\oorank(X,E)$ and the following classes of mappings: $\Bof_1(X,E)$, $\Hf_1(X,E)$, $\Frag(X,E)$, mappings with the PCP and mappings of Baire class $1$. The positive results of this kind are contained in Theorem~\ref{porovnani}.

\begin{lemma}\label{omegabof}
Let $X$ be a topological space, $E$ a metric space and $f:X\to E$ with $\orank(f)<\omega_1$. Then $f\in \Bof_1(X,E)$.
\end{lemma}
\begin{proof}
For every $n\in\en$ we set $\alpha_n = \orank(f,\frac1n)<\omega_1$ and find the open sets $\{ U_\alpha^n:\alpha\leq\alpha_n\}$ from the definition of the oscillation rank. Then, we define the families
\[
\mathcal{D}_n^{\alpha} = \{ G\subset X: G \mbox{ is open, } \diam f(G\setminus U_\alpha^n)<\frac1n\}, \quad n\in\en, \alpha\leq\alpha_n,
\]
and set $\mathcal{D} = \bigcup_{n\in\en} \mathcal{D}_n^{\alpha}$.

Now, let $U\subset E$ be an open set. Then clearly
\[
f^{-1}(U) = \bigcup_{n\in\en}\bigcup_{\alpha\leq\alpha_n} \bigcup_{G\in \mathcal{D}_n^{\alpha}} \left\{ G\setminus U_\alpha^n : f(G\setminus U_\alpha^n)\subset U\right\} = \bigcup_{n\in\en}\bigcup_{\alpha\leq\alpha_n} G_n^{\alpha}\setminus U_\alpha^n,
\]
where $G_\alpha^n = \bigcup\{ G\in\mathcal{D}_n^{\alpha}: f(G\setminus U_\alpha^n)\subset U\}$ for every $n\in\en$ and $\alpha\leq\alpha_n$.

Since $\alpha_n<\omega_1$ for every $n\in\en$, it follows that $f\in\Bof_1(X,E)$.
\end{proof}

\begin{thm}\label{porovnani}
Let $X$ be a topological space and $E$ a metric space. Then the following schema holds true
\[
\begin{matrix}
\Hf_1(X,E) & \supset & \Bof_1(X,E)\\
\cup & & \cup\\
\Frag(X,E) & \supset & \oorank(X,E)
\end{matrix}
\]
and no other nontrivial inclusion is valid generally.
\end{thm}

\begin{proof}
Inclusion $\Frag(X,E)\subset\Hf_1(X,E)$ was established by \cite[Theorem 2.3]{koumoullis}. 

Inclusion $\Bof_1(X,E)\subset\Hf_1(X,E)$ is a direct consequence of definitions. Inclusion $\oorank(X,E)\subset\Frag(X,E)$ follows from Lemma~\ref{nutnost} and inclusion $\oorank(X,E)\subset\Bof_1(X,E)$ is exactly Lemma~\ref{omegabof}. 

Assuming Martin's Axiom and the negation of the continuum hypothesis there is a function from $\Bof_1(X,E)$ constructed in \cite[Example 2.4(3)]{koumoullis} which is not fragmented. Therefore $\Bof_1(X,E)\nsubseteq\Frag(X,E)$. In \cite[Remark 3.3]{spurny-amh} there is a space constructed containing a resolvable non-Borel set $A$. The characteristic function $\chi_A$ of the set $A$ is clearly fragmented and $\chi_A\notin \Bof_1(X,E)$. Hence, $\Frag(X,E)\nsubseteq\Bof_1(X,E)$.

\end{proof}


From the viewpoint of Theorem~\ref{porovnani} and \cite[Theorem 2.3]{koumoullis} it seems reasonable to investigate the relation between the mappings with the PCP and the mappings with countable oscillation rank.

\begin{example}
There exists a function $f:\qe\to\er$ lacking the PCP such that $f\in \oorank(\qe,\er)$.
\end{example}
\begin{proof}
The function $f$ from \cite[Examples 2.4(4)]{koumoullis} does the job.
\end{proof}

\begin{example}
There exists a topological space $X$, a metric space $E$ and a mapping $f:X\to E$ with the PCP and satisfying $\orank(f)= \omega_1$.
\end{example}
\begin{proof}
Lemma~\ref{prostory} provides a hereditarily Baire topological space $X$ and a function $f:X\to \{0,1\}$ such that $\orank(f) =\omega_1$ and $f\in\Bof_1(X,\{0,1\})$. 

Then $f$ is fragmented by Lemma~\ref{nutnost} (and hence is a characteristic function of a resolvable set). Applying \cite[Theorem 2.3]{koumoullis}, $f$ has the PCP.
\end{proof}

Let us mention that if $f$ is the characteristic function of a resolvable non-Borel set in a topological space (see \cite[Remark 3.3]{spurny-amh}), then $f$ has the PCP and by Lemma~\ref{omegabof} satisfies $\orank(f)\geq\omega_1$.

\medskip

Let $K$ be a compact metrizable space. Then, according to \cite[Proposition 2]{kechris1990classification}, a function $f:K\to\er$ is of Baire class $1$ if and only if $\beta(f)<\omega_1$. A question arises whether there is a chance for an analogous proposition to hold in a more general setting. Alas, there is no relation between functions of Baire class $1$ and the class $\oorank(X,\er)$ for a general topological space $X$.

\begin{example}
There exists a function $g:\qe\to\er$ of Baire class $1$ which is not fragmented.
\end{example}
\begin{proof}
The function $g$ from \cite[Examples 2.4(4)]{koumoullis} has the desired properties.
\end{proof}

\begin{example}
There exists a topological space $X$ and a function $f:X\to\er$ such that $f\in\oorank(X,\er)$ but which is not of Baire class $1$.
\end{example}
\begin{proof}
Let $X$ be a topological space containing an open set which is not $F_\sigma$. Then the characteristic function $f$ of this set has the desired properties.
\end{proof}

\medskip

\noindent\textbf{Aknowledgements.} The author would like express his thanks to Ji{\v{r}}{\'{i}} Spurn{\'{y}} and Ond\v{r}ej Kalenda for introducing him the topic and many useful suggestions.

The author were supported by the Grant IGA\_PrF\_2022\_008 ``Mathematical Models'' of the Internal Grant Agency of Palack\'y University in Olomouc.

\bibliography{Ludvik-submission}\bibliographystyle{siam}

\begin{thebibliography}{10}

\bibitem{angosto-dis}
{\sc C.~Angosto}, {\em Distancia a espacios de funciones}, PhD thesis,
  Universidad de Murcia, Spain, 2007.

\bibitem{quantitative}
{\sc C.~Angosto and B.~Cascales}, {\em The quantitative difference between
  countable compactness and compactness}, J. Math. Anal. Appl., 343 (2008),
  pp.~479--491.

\bibitem{ca-b1}
{\sc C.~Angosto, B.~Cascales, and I.~Namioka}, {\em Distances to spaces of
  {B}aire one functions}, Math. Z., 263 (2009), pp.~103--124.

\bibitem{ca-ma-ra}
{\sc B.~Cascales, W.~Marciszewski, and M.~Raja}, {\em Distance to spaces of
  continuous functions}, Topology Appl., 153 (2006), pp.~2303--2319.

\bibitem{engelking}
{\sc R.~Engelking}, {\em General topology}, vol.~6 of Sigma Series in Pure
  Mathematics, Heldermann Verlag, Berlin, second~ed., 1989.
\newblock Translated from the Polish by the author.

\bibitem{fabian2010banach}
{\sc M.~Fabian, P.~Habala, P.~H{\'a}jek, V.~Montesinos, and V.~Zizler}, {\em
  Banach space theory, The basis for linear and nonlinear analysis}, CMS Books
  in Mathematics/Ouvrages de Math\'ematiques de la SMC, Springer, New York,
  2011.

\bibitem{Ha}
{\sc R.~W. Hansell}, {\em Borel measurable mappings for nonseparable metric
  spaces}, Trans. Amer. Math. Soc., 161 (1971), pp.~145--169.

\bibitem{hansell1974characterizing}
\leavevmode\vrule height 2pt depth -1.6pt width 23pt, {\em On characterizing
  non-separable analytic and extended {B}orel sets as types of continuous
  images}, Proc. London Math. Soc. (3), 28 (1974), pp.~683--699.

\bibitem{hansell1998nonseparable}
{\sc R.~W. Hansell}, {\em Nonseparable analytic metric spaces and quotient
  maps}, Topology Appl., 85 (1998), pp.~143--152.
\newblock 8th Prague Topological Symposium on General Topology and Its
  Relations to Modern Analysis and Algebra (1996).

\bibitem{Haydon91oncertain}
{\sc R.~Haydon, E.~Odell, and H.~Rosenthal}, {\em On certain classes of
  {B}aire-{$1$} functions with applications to {B}anach space theory}, in
  Functional analysis ({A}ustin, {TX}, 1987/1989), vol.~1470 of Lecture Notes
  in Math., Springer, Berlin, 1991, pp.~1--35.

\bibitem{jaynefrag}
{\sc J.~E. Jayne, J.~Orihuela, A.~J. Pallar{\'e}s, and G.~Vera}, {\em
  {$\sigma$}-fragmentability of multivalued maps and selection theorems}, J.
  Funct. Anal., 117 (1993), pp.~243--273.

\bibitem{K97}
{\sc O.~Kalenda}, {\em Note on connections of the point of continuity property
  and {K}uratowski problem on function having the {B}aire property}, Acta Univ.
  Carolin. Math. Phys., 38 (1997), pp.~3--12.

\bibitem{K98}
{\sc O.~Kalenda}, {\em New examples of hereditarily {$t$}-{B}aire spaces},
  Bull. Polish Acad. Sci. Math., 46 (1998), pp.~197--210.

\bibitem{kechris1990classification}
{\sc A.~S. Kechris and A.~Louveau}, {\em A classification of {B}aire class
  {$1$} functions}, Trans. Amer. Math. Soc., 318 (1990), pp.~209--236.

\bibitem{koumoullis}
{\sc G.~Koumoullis}, {\em A generalization of functions of the first class},
  Topology Appl., 50 (1993), pp.~217--239.

\bibitem{kuratowski-top}
{\sc K.~Kuratowski}, {\em Topology. {V}ol. {I}}, New edition, revised and
  augmented. Translated from the French by J. Jaworowski, Academic Press, New
  York, 1966.

\bibitem{kuratowski-set}
{\sc K.~Kuratowski and A.~Mostowski}, {\em Set theory}, Translated from the
  Polish by M. Maczy\'nski, PWN-Polish Scientific Publishers, Warsaw, 1968.

\bibitem{fine}
{\sc J.~Luke{\v{s}}, J.~Mal{\'y}, and L.~Zaj{\'{\i}}{\v{c}}ek}, {\em Fine
  topology methods in real analysis and potential theory}, vol.~1189 of Lecture
  Notes in Mathematics, Springer-Verlag, Berlin, 1986.

\bibitem{marsol}
{\sc D.~A. Martin and R.~M. Solovay}, {\em Internal {C}ohen extensions}, Ann.
  Math. Logic, 2 (1970), pp.~143--178.

\bibitem{Ro74}
{\sc H.~P. Rosenthal}, {\em A characterization of {B}anach spaces containing
  {$l\sp{1}$}}, Proc. Nat. Acad. Sci. U.S.A., 71 (1974), pp.~2411--2413.

\bibitem{lattices}
{\sc H.~H. Schaefer}, {\em Banach lattices and positive operators},
  Springer-Verlag, New York, 1974.
\newblock Die Grundlehren der mathematischen Wissenschaften, Band 215.

\bibitem{Scott61}
{\sc D.~Scott}, {\em Measurable cardinals and constructible sets}, Bull. Acad.
  Polon. Sci. S\'{e}r. Sci. Math. Astronom. Phys., 9 (1961), pp.~521--524.

\bibitem{Shoenfield59}
{\sc J.~R. Shoenfield}, {\em On the independence of the axiom of
  constructibility}, Amer. J. Math., 81 (1959), pp.~537--540.

\bibitem{spurny-amh}
{\sc J.~Spurn{\'y}}, {\em Borel sets and functions in topological spaces}, Acta
  Math. Hungar., 129 (2010), pp.~47--69.

\bibitem{spurny-studia}
\leavevmode\vrule height 2pt depth -1.6pt width 23pt, {\em Distances to spaces
  of affine {B}aire-one functions}, Studia Math., 199 (2010), pp.~23--41.

\end{thebibliography}

\end{document}